\documentclass{conm-p-l}
\usepackage{amscd, amssymb}
 
\theoremstyle{plain}

\newtheorem{lemma}{Lemma}
\newtheorem*{mainlemma}{Main Lemma} 
\newtheorem{theorem}{Theorem}
\newtheorem*{main}{Main Theorem}
\newtheorem{corollary}{Corollary}
\newtheorem*{claim*}{Claim}

\newtheorem{question}{Question}

\theoremstyle{definition}

\newtheorem*{definition*}{Definition}

\theoremstyle{remark}
\newtheorem{remark}{Remark}
\newtheorem*{remark*}{Remark}

\newtheorem*{acknowledgement}{Acknowledgement}

\newcommand{\cc}{\mathbb{C}}

\newcommand{\mbar}{{\overline{M}}}

\newcommand{\p}{\partial}

\newcommand{\h}{\hbar}

\newcommand{\bt}{t}

\newcommand{\HH}{\mathcal{H}}

\newcommand{\on}{\operatorname}

\renewcommand{\L}{\mathcal{L}}
\newcommand{\bq}{q}
\newcommand{\bp}{p}

\newcommand{\f}{f}
\newcommand{\g}{g}
\newcommand{\1}{\mathbf{1}}

\newcommand{\rspin}{\text{$r$-spin}}
\newcommand{\WK}{\text{WK}}

\title[Witten's conjecture and Virasoro conjecture]
{Witten's conjecture and Virasoro conjecture\\
for genus up to two}

\author{Y.-P.~Lee} 

\address{Department of Mathematics \\
        University of Utah \\
        Salt Lake City, Utah 84112-0090}
\email{yplee@math.utah.edu}

\thanks{Research partially supported by NSF.}

\begin{document}

\begin{abstract}
This is an expository paper based on the results in \cite{YL2} and \cite{GL}.
The main goal is to show how one may prove the following two conjectures
\emph{for genus up to two}:
\begin{enumerate}
\item Witten's conjecture on the relations between higher spin curves and 
Gelfand--Dickey hierarchy.
\item Virasoro conjecture for conformal semisimple Frobenius manifolds. 
\end{enumerate}
\end{abstract}

\maketitle

\section{Introduction}

\subsection{The conjectures}

\subsubsection{Witten's conjecture}

E.~Witten in 1990 made a striking conjecture between generating functions
of intersection numbers on moduli spaces of stable curves and
a $\tau$-function of KdV hierarchies \cite{EW1}. This conjecture
says that the following geometrically defined function
\[
 \tau^{pt}(t_0,t_1, \ldots) 
 = e^{\sum_{g=0}^{\infty} \hbar^{g-1} F^{pt}_g(t_0,t_1, \ldots)}
\]
is a $\tau$-function of the KdV hierarchy.
\footnote{$\hbar$ is usually put to be 1.} 
In the above formula, 
$F^{pt}_g(t)$ is the generating function of (tautological) intersection
numbers on the moduli space of stable curves of genus $g$.
Moreover, from elementary geometry
of moduli spaces, one easily deduces that $\tau^{pt}$ satisfies an additional 
equation, called the \emph{string equation}. 
It is known from the theory of KdV hierarchy that the string equation 
for the KdV (or in general KP) hierarchies uniquely determines a 
$\tau$-function parameterized by the Sato's grassmannian.
This particular $\tau$-function will be called Witten--Kontsevich 
$\tau$-function and denoted $\tau_{\WK}$. 
In other words, $\tau_{\WK}=\tau^{pt}$.
Often $\tau^{pt}$ is used to emphasize its geometric nature and $\tau_{\WK}$
is used when the integrable system side is emphasized. 

In 1991 Witten formulated a remarkable generalization of the above conjecture.
He argued that an analogous generating function $\tau^{\rspin}$
of the intersection numbers on moduli spaces of $r$-spin curves should be 
identified as a $\tau$-function of $r$-th Gelfand--Dickey ($r$-KdV) hierarchies
\cite{EW2}. When $r=2$, this conjecture will \emph{eventually} reduces to the 
previous one (after some work), as $2$-KdV is the ordinary KdV.

The special case was soon proved by M.~Kontsevich \cite{MK}. 
More recently a new proof was given by Okounkov--Pandharipande \cite{OP}. 
However, the generalized conjecture remains open up to this day.

It may be worth pointing out that the status of two conjectures was very 
different when they were first proposed. 
The 1990 conjecture was from the beginning formulated mathematically, 
using only well defined mathematical quantities. The 1991 conjecture,
on the contrary, involves the concepts like moduli spaces of $r$-spin
curves and their virtual fundamental classes (in modern terminologies) 
for which Witten offers only sketches of their construction.
Perhaps the sharpest contrast lies in the fact that
there were plenty of evidences supporting the 1990 conjecture,
but virtually no evidences supporting 1991 conjecture beyond genus zero at the 
time they are formulated. 

Throughout the years, T.~Jarvis, and later joint by T.~Kimura, A.~Vaintrob,
and A.~Polishchuk, T.~Mochizuki, have clarified the foundational issues.
In particular, Jarvis--Kimura--Vaintrob \cite{JKV} established the genus zero 
case of the conjecture; Mochizuki and Polishchuk
independently established the following property for $\tau^{\rspin}$:

\begin{theorem} \cite{TM} \cite{AP}
All tautological equations hold for $F_g^{\rspin}$. 
\end{theorem}
In fact, $F_g^{\rspin}$ satisfies all ``expected functorial properties'',
similar to the axioms formulated by Kontsevich--Manin in the Gromov--Witten
theory. 
Note that he term \emph{tautological equations} usually mean the relations of 
tautological (Chow) classes on moduli spaces of curves. 
Here they are used for the relations of the classes on the moduli spaces of
spin curves, obtained by the forgetful map $\mbar_{g,n}^{1/r} \to \mbar_{g,n}$.
Later, they are also used for the induced relations in the Gromov--Witten 
theory.

However, Riemann's trichotomy of Riemann surfaces has taught us that things
are very different in genus one and at higher genus.
Our Main Theorem therefore provides a solid confirmation for Witten's 1991
conjecture, covering one example ($g=1$ and $g=2$) for the other two cases
in the trichotomy. In fact, this work starts as a project trying to understand 
this conjecture at higher genus.

For more background information about Witten's conjecture, the readers are 
referred to Witten's original article \cite{EW2} and the paper \cite{JKV} by 
Jarvis--Kimura--Vaintrob, both well-written.
In the remaining of this article, ``Witten's conjecture'' means the 1991 
conjecture if not otherwise specified.

\subsubsection{Virasoro conjecture}

In 1997 another generalization of Witten's 1990 conjecture was proposed by
T.~Eguchi, K.~Hori and C.~Xiong \cite{EHX}. Witten's 1990 conjecture has 
an equivalent formulation: $\tau^{pt}$ is annihilated by infinitely many 
differential operators $\{ L^{pt}_n \}, n \ge -1$, satisfying the Virasoro
relations
\[
  [L_m, L_n]=(m-n) L_{m+n}. 
\]
Eguchi--Hori--Xiong, and S.~Katz, managed to formulate a conjecture for
any projective smooth variety $X$, generalizing the above assertion. 
Namely, they found the formulas of $\{ L^X_n \}$ for $n \ge -1$, satisfying 
Virasoro relations and conjectured that
\[
  L^X_n \tau^X (\bt) =0, \quad \text{for $n \ge -1$}.
\]
In the above formula,
\[
 \tau^X (\bt) := e^{\sum_{g=0}^{\infty} \hbar^{g-1} F^X_g(\bt)},
\]
and $F^X_g (\bt)$ is the generating function of genus $g$ Gromov--Witten
invariants with descendents for the projective manifold $X$.
This conjecture is commonly referred to as the Virasoro conjecture.

Eguchi--Hori--Xiong was able to give strong evidences for their conjecture at
genus zero. Later X.~Liu and G.~Tian \cite{LT} established the genus zero case.
Using a very different method, B.~Dubrovin and Y.~Zhang established the
genus one case of Virasoro conjecture for 
\emph{conformal} semisimple Frobenius manifolds.
\footnote{The definition of Frobenius manifolds in this article does not
require existence of an Euler field, which is assumed in Dubrovin's definition.
Dubrovin's definition will be referred to as 
\emph{conformal Frobenius manifold} instead.}

The recent major developments are Givental's proof of Virasoro conjecture for
toric Fano manifolds \cite{AG2} and Okounkov--Pandharipande's proof of Virasoro
conjecture for algebraic curves \cite{OP2}. 

\subsection{Givental's theory}

A.~Givental introduces a completely new approach to Gromov--Witten theory
in a series of papers \cite{AG1, AG2, AG4}, dating back to August 2000.
It is beyond our ability to summarize Givental's theory in a few paragraphs. 
In the following we will restrict ourselves to some highlights of his theory,
mainly for the purpose of fixing the notations and putting the theorems in
a framework. 
The details can be found in \cite{LP}.

The essence of his theory is a construction of a ``combinatorial model'' of
higher genus invariants via graphic enumeration, with the information of
edges coming from the underlying semisimple Frobenius manifold 
(i.e.~genus zero theory) and information of vertices from $\tau^{pt}$.
Formulaically, given a semisimple Frobenius manifold $H$ of dimension $N$,
he defines an operator $\hat{O}_H = \exp (\hat{o}_H)$.
Givental's $\tau$-function is defined to be
\begin{equation} \label{e:taug}
 \tau_G^H: = e^{\sum_{g=0}^{\infty} \hbar^{g-1} G^H_g (t^1,\ldots, t^N)} 
 := \hat{O}_H \prod_{i=1}^N \tau_{\WK}(t^i, \hbar) .
\end{equation}
In fact, the operator $\hat{O}_H$ is a special kind of operator and belongs to
quantized twisted loop group, which will be discussed in Section~\ref{s:3}.
\footnote{$\hat{O}_H$ is actually not really in the quantized twisted loop 
group, but in its completion. We will ignore the difference in this article.}
The Feynman rules then dictate a formula for $G_g$. When the Frobenius
manifold comes from geometry, i.e.~$H=QH^*(X)$, Givental conjectures that
his combinatorial model is the same as the geometric model. That is
$G^H_g=F^X_g$ when $H=QH^*(X)$.

What makes Givental's model especially attractive are the facts that 
\begin{enumerate}
\item it works for any semisimple Frobenius manifolds
\item it enjoys properties often complementary to the geometric theory. 
\end{enumerate}
Thanks to (1), one also has a Givental's model for the Frobenius manifolds
$H_{A_{r-1}}$ of the miniversal deformation space of $A_{r-1}$ singularity. 
It turns out that this Frobenius manifold is isomorphic to the Frobenius 
manifold defined by the genus zero potential of $r$-spin curves. 
Furthermore, Givental has recently proved

\begin{theorem} \cite{AG3} \label{t:ag3}
$\tau_G^{H_{A_{r-1}}}$ is a $\tau$-function of $r$-KdV hierarchy.
\end{theorem}

As in the case of the ordinary KdV, it is easy to show that both
$\tau_G^{H_{A_{r-1}}}$ and $\tau^{\rspin}$ satisfy the additional string 
equation. Therefore, in order to prove Witten's conjecture, one only has
to answer the following question positively,

\begin{question} \label{q:1}
Is $G_g^{H_{A_{r-1}}}= F^{\rspin}_g$?
\end{question}

As examples for (2), $\tau_G^H$ satisfies Virasoro conjecture 
and $\tau^X$ satisfies the tautological relations 
almost by definitions. In the second case, one notes that if the theory is
defined geometrically, one can pull-back the relations on moduli spaces of
curves to moduli spaces of maps. In the first case, one defines the Virasoro 
operators for semisimple Frobenius manifold $H$ by 
\[
 L_n^H (\bt) := \hat{O} \prod_i L^{pt}_n(t_i) \hat{O}^{-1},
\]
it is obvious that $L_n^H$ also satisfy Virasoro relations.
One immediately gets Virasoro conjecture for $H$ by Kontsevich's theorem.
It is also true that $L_n^H = L_n^X$ when the semisimple Frobenius manifold
$H$ comes from quantum cohomology of $X$, i.e.~$H=QH^*(X)$.

However, the converse statements pose nontrivial challenges. 

\begin{question} 
Does the tautological relations hold for $G_g$?
\end{question}

\begin{question}
Does Virasoro conjecture hold for $\tau^X$?
\end{question}

An obvious, and indeed very good, strategy to answer all the above three
questions is to answer a generalized version of Question~1:

\begin{question} \label{q:4}
Is $G_g = F_g$? That is, does the combinatorial construction coincide with
the geometric one when both are available?
\end{question}

A positive answer to Question~4 obviously answer all three questions at once.
Interesting enough, the solution to Question~4 turns out to be closely 
related to, an in many cases equivalent to, the solutions of Questions~2 and 3
by some uniqueness theorems.

For example, the answer to Question~3 is equivalent to that to Question~4 in 
the geometric Gromov--Witten theory by a result of Dubrovin and Zhang 
\cite{DZ3}. They proved that Virasoro conjecture plus $(3g-2)$-jet conjecture 
(actually a theorem of E.~Getzler \cite{EG2} in geometric Gromov--Witten theory
and of Givental in the context of semisimple Frobenius manifolds) uniquely 
determines $\tau$-function for any semisimple Frobenius manifold.
It is expected that $\tau^{\rspin}$ will also satisfy the $(3g-2)$-jet 
property. 

As for Question~2, the equivalence is established by some uniqueness theorems
in genus one and two. Similar phenomena are ``expected'' to hold in higher 
genus as well, although there is no hard evidences at this moment.
In genus one, Dubrovin and Zhang made the following important observation
about the uniqueness theorem. 

\begin{lemma} \cite{DZ1} \label{l:dz}
The genus one descendent potentials for any semisimple Frobenius manifolds $H$
are uniquely determined, up to linear combination of canonical coordinates, 
by genus zero potentials, genus one topological recursion relations, 
and genus one Getzler's equation.

Furthermore, if $H$ is conformal, then the genus one potential is uniquely
determined up to constant terms.
\end{lemma}

The proof of this fact goes as follows. 
First, genus one TRR guarantees that the descendent invariants are uniquely 
determined by primary invariants. 
Second, genus one Getzler's equation, when written in canonical coordinates 
$u^i$, is equal to $\frac{\p^2 F_1}{\p u^i \p u^j} = B_{ij}$ where $B_{ij}$ 
involves only genus zero invariants. 
Moreover, the conformal structure determined by
a linear vector field (Euler field), uniquely determines the linear term.

The uniqueness theorem in genus two, proved by X.~Liu, is much more involved:

\begin{theorem} \cite{XL} \label{t:xl}
The genus two descendent potentials for any conformal semisimple Frobenius 
manifolds are uniquely determined up to constants by genus two equations by
Mumford, Getzler \cite{EG} and Belorousski--Pandharipande (BP) \cite{BP}.
\end{theorem}

It is worth noting that whether this uniqueness theorem, or any weaker version,
holds for non-conformal semisimple Frobenius manifolds remains unknown.

\subsection{Statements of the main results}

By Lemma~\ref{l:dz}, the differentials of the genus one potentials are 
equal, $d G_1=d F_1$, if $G_1$ satisfies genus one TRR and genus one Getzler's
equation, plus some initial condition to fix the constant terms. The positive
answer to Question~2 in genus one and therefore all other questions are proved 
in \cite{GL}

\begin{theorem} \cite{GL} \label{t:g1}
$d G_1 = d F_1$ for all semisimple Frobenius manifolds.
\end{theorem}

This theorem generalizes the earlier results by Dubrovin--Zhang for
\emph{conformal} semisimple Frobenius manifolds in \cite{DZ1}.

\begin{theorem} \cite{YL2} \label{t:g2}
$G_2$ satisfies genus two tautological equations by Mumford, Getzler and BP.
\end{theorem}

These two theorems, combined with the above results, immediately implies 

\begin{main} \cite{YL2}
Witten's conjecture and
Virasoro conjecture for conformal semisimple Frobenius manifolds
hold up to genus two.
\end{main}

\subsubsection{Main ideas involved in the proofs}

To prove Theorem~\ref{t:g1} and \ref{t:g2}, note that 
\begin{itemize}
\item $\tau_G = \hat{O} \prod_i \tau^{pt}(\bt^i)$.
\item $\tau^{pt}(\bt)$ satisfies all tautological equations.
\end{itemize}

Therefore, in order to prove $G_g$ satisfies tautological equations, one only
has to prove that these equations are invariant under the action of $\hat{O}$.
This is the approach taken in \cite{GL} and \cite{YL2}.

\begin{remark}
There are other possible approaches to this problem.
Our earlier approach in \cite{YL1} reduces the checking of 
Theorem~\ref{t:g2} to a complicated, but finite-time checkable, identities.
Nevertheless, it lacks the underlying simplicity of this approach.

After this result was announced, X.~Liu informed us (and later posted in arxiv
\cite{XL2}) that he was also able to prove the genus two Virasoro conjecture 
for the conformal semisimple Gromov--Witten theory by reducing it
to some complicated identities which he was able to check by hand \emph{and}
by a Mathematica programs.
\end{remark}

\begin{acknowledgement}
I would like to thank B.~Dubrovin, T.~Jarvis, T.~Kimura, X.~Liu, 
Y.~Ruan and A.~Vaintrob for useful discussions and communications. I am also
thankful to E.~Getzler for pointing out an inaccurate statement in an earlier 
version. 
Many ideas of this work comes from my collaborations with A.~Givental
\cite{GL} and R.~Pandharipande \cite{LP}. It is a great pleasure to thank 
both of them.
\end{acknowledgement}

\section{Frobenius manifolds}

\subsection{Givental's theory of formal Frobenius manifolds}

Let $H$ be a complex vector space of dimension $N$ with a distinguished element
$\1$. Let $(\cdot,\cdot)$ be a $\cc$-bilinear metric on $H$, 
i.e.~a nondegenerate symmetric $\cc$-bilinear form.
Let $\HH$ denote the (infinite dimensional) complex vector space $H((z^{-1}))$ 
consisting of Laurent formal series in $1/z$ with vector coefficients.
\footnote{Different completions of this spaces are used in different places,
but this subtlety will be ignored in the present article.}
Introduce the symplectic form $\Omega$ on $\HH$:
\[ \Omega (\f,\g )  = \frac{1}{2\pi i} \oint\ ( \f (-z), \g(z) ) \ dz .\]
The polarization $\HH = \HH_{+}\oplus \HH_{-}$ by the Lagrangian subspaces
$\HH_{+} = H[z]$ and $\HH_{-} = z^{-1} H [[ z^{-1} ]]$ provides a symplectic
identification of $(\HH, \Omega)$ with the cotangent bundle $T^*\HH_{+}$.

Let $\{ \phi_{\mu} \}$ be an \emph{orthonormal} basis of $H$. 
Introduce Darboux coordinates $\{ \bp_k^{\mu}, \bq_k^{\mu} \}$, 
$k= 0, 1, 2, \ldots$ and $\mu=1, \ldots, N$, compatible with this polarization,
so that
\[
 \Omega = \sum d \bp^{\mu}_k \wedge d \bq^{\mu}_k .
\] 
An $H$-valued Laurent formal series can be written in this basis as
\begin{multline*}
 \ldots + (\bp^1_1,\ldots,\bp^N_{1}) \frac{1}{(-z)^2}
 + (\bp^1_{0},\ldots, \bp^N_{0}) \frac{1}{(-z)} \\
 + (\bq^1_{0},\ldots, \bq^N_{0})
 + (\bq^1_{1}, \ldots, \bq^N_{1}) z + \ldots.
\end{multline*}
To simplify the notations, $\bp_k$ will stand for the vector
$(\bp^1_{k},\ldots,\bp^N_{k})$ and $\bp^{\mu}$ for 
$(\bp^{\mu}_0, \bp^{\mu}_1, \ldots )$.

Let $A(z)$ be an $\on{End}(H)$-valued Laurent formal series in $z$ satisfying
\[
  (A(-z) f(-z), g(z)) + (f(-z), A(z) g(z)) =0,
\]
then $A(z)$ defines an infinitesimal symplectic transformation
\[
  \Omega(A f, g) + \Omega(f, A g)=0.
\]
An infinitesimal symplectic transformation $A$ of $\HH$ corresponds to a
quadratic polynomial $P(A)$ in $\bp, \bq$
\[
  P(A)(f) := \frac{1}{2} \Omega(Af, f) .
\]
($A$ is a symplectic vector fields on the symplectic vector space $(\HH,
\Omega)$, and the relation between the function $P(A)$ and vector field $A$ is
$d P(A) = i_A \Omega$.) For example, if $\dim H=1$ and $A(z) =1/z$, then
\begin{equation*} 
 P(z^{-1})= -\frac{\bq_0^2}{2} - \sum_{m=0}^{\infty} \bq_{m+1} \bp_m.
\end{equation*}



\subsection{Lagrangian cones}

Let $F_0(\bt)$ be a formal series in $\bt$, where
$\bt = (\bt_0,\bt_1,\bt_2,...)$ is related to 
$\bq = (\bq_0, \bq_1, \bq_2,\ldots)$ through the following change of variables:
\[
 \sum_{k=0}^{\infty} \bq_k z^k =: \bq (z)
 = \bt (z) -z \1 := - z \1 + \sum_{k=0}^{\infty} \bt_k z^k .
\]
Thus the formal function $F_0 (\bt(z))$ near $\bt=0$ 
becomes a formal function $F_0(\bq)$ on the space $\HH_{+}$ near the point
$\bq (z)= -z$. This convention is called the \emph{dilaton shift}.

In the Gromov--Witten theory, $F_0(\bt)$ is the genus zero descendent 
potential. It satisfies many properties due to the geometry of the moduli 
spaces. Three classes of partial differential equations are most relevant.
They are called  the \emph{Topological Recursion Relations} (TRR), 
the {\em String Equation} (SE) and the {\em Dilaton Equation} (DE):
\begin{align}
   \tag{DE} 
 &\frac{\p F_0(\bt)}{\p \bt^1_1} (\bt)  
 = \sum_{n=0}^{\infty} \sum_{\nu} \bt^{\nu}_n 
 \frac{\p F_0(\bt)}{\p \bt^{\nu}_n} - 2F_0 (\bt), \\
 \tag{SE}  
  &\frac{\p F_0 (\bt)}{\p \bt^1_0} =
  \frac{1}{2}(\bt_0,\bt_0)+ \sum_{n=0}^{\infty}
  \sum_{\nu} \bt_{n+1}^{\nu} \frac{\p F_0 (\bt)}{\p \bt^{\nu}_n}, \\
  \label{TRR} \tag{TRR}
  &\frac{\p^3 F_0 (\bt)}
       {\p \bt^{\alpha}_{k+1} \p \bt^{\beta}_{l} \p \bt^{\gamma}_m}
  = \sum_{\mu} \frac{\p^2 F_0 (\bt)}{\p \bt^{\alpha}_{k} \p \bt^{\mu}_{0}}
    \frac{\p^3 F_0 (\bt)}
       {\p \bt^{\mu}_{0} \p \bt^{\beta}_{l} \p \bt^{\gamma}_m}
\end{align}
for all $\alpha,\beta,\gamma$ and all $k,l,m\geq 0$.
Note that $t^1_k$ denote the dual coordinates of the vectors $\1 z^k$.

Denote by $\L$ the graph of the differential $dF_0$:
\[ 
 \L = \{ (\bp, \bq) \in T^*\HH_{+}:\ \bp = d_{\bq} F_0 \} .
\]
It is considered as a formal germ at $\bq = -z$ (i.e.~$\bt=0$) of a Lagrangian
section of the cotangent bundle $T^*\HH_{+}$ and can therefore be considered as
a formal germ of a Lagrangian submanifold in the symplectic loop space
$(\HH,\Omega)$.

\begin{theorem} \cite{AG4} \label{t:1}
The function $F_0$ satisfies TRR, SE and DE if and only if the corresponding
Lagrangian submanifold $\L \subset \HH$ has the following properties:
\begin{enumerate}
\item $\L$ is a Lagrangian cone with the vertex at the origin.
\item The tangent spaces $L_{\f}=T_{\f}\L$ satisfy $zL_{\f} \subset L_{\f}$ 
(and therefore $\dim L_{\f}/zL_{\f} = \dim \HH_{+} /z\HH_{+} =\dim H$).
\item $zL_{\f} \subset \L$.
\item The same $L_{\f}$ is the tangent space to $\L$ not only along the line 
of $\f$ but also at all smooth points in $zL_{\f} \subset \L$. 
\end{enumerate}
\end{theorem}


One may rephrase the above properties by saying that $\L$ is a cone ruled by
the isotropic subspaces $zL$ varying in a $\dim H$-parametric family with the 
tangent spaces along $zL$ equal to the same Lagrangian space $L$. 
This in particular implies that the family of $L$ generates a 
{\em variation of semi-infinite Hodge structures} 
in the sense of S.~Barannikov, i.e. a family of semi-infinite flags
\[  \cdots zL \subset L \subset z^{-1}L \cdots \]
satisfying the Griffiths integrability condition.

We note the following theorem 

\begin{theorem} \cite{AG4}
A Lagrangian cone satisfying the above conditions defines a formal germ of a
Frobenius manifold plus a given calibration, in the sense of \cite{AG2}, and
vice versa. 
\end{theorem}

Although two formulations are equivalent, the Lagrangian cone formulation
is much more transparent and geometric. One can say that the Lagrangian
cone formulation is the geometrization of the equations SE, DE, and TRR.
Moreover, these properties are formulated in terms of the symplectic structure 
$\Omega$ and the operator of multiplication by $z$. 
Hence it does not depend on the choice of the polarization. 
This shows that the system DE$+$SE$+$TRR has a huge symmetry group.

\begin{definition*}
Let $\L^{(2)}GL(H)$ denote the {\em twisted loop group} which consists of 
$\operatorname{End}(H)$-valued formal Laurent formal series $M(z)$ in
the indeterminate $z^{-1}$ satisfying $M^*(-z)M(z)=\1$. 
Here $\ ^*$ denotes the adjoint with respect to $(\cdot ,\cdot )$.
\end{definition*} 

The condition $M^*(-z)M(z)=\1$ means that $M(z)$ is a symplectic transformation
on $\HH$.

\begin{corollary} \label{c:tlg}
The action of the twisted loop group preserves the class of the Lagrangian
cones $\L$ satisfying (1-4) and, generally speaking, yields new generating
functions $F_0$ which satisfy the system DE$+$SE$+$TRR whenever the original 
one does.
\end{corollary}

\section{Higher genus and quantization} \label{s:3}

To quantize an infinitesimal symplectic transformation, or its corresponding
quadratic hamiltonians, we recall the standard Weyl quantization. 
A polarization $\HH=T^* \HH_+$ on the symplectic vector space $\HH$ 
(the phase space) defines a configuration space $\HH_+$.
The quantum ``Fock space'' will be a certain class of functions $f(\h, q)$ on 
$\HH_+$ (containing at least polynomial functions), with additional formal
variable $\hbar$ (``Planck's constant'').
The classical observables are certain functions of $\bp, \bq$.
The quantization process is to find for the classical mechanical system on 
$\HH$ a ``quantum mechanical'' system on the Fock space such that
the classical observables become operators on the Fock space. 
In particular, the classical hamiltonians $h(q,p)$ on $\HH$ are quantized
to be differential operators $\hat{h}(q,\dfrac{\p}{\p q})$ on the Fock space. 

In the above Darboux coordinates, the quantization $P \mapsto \hat{P}$ assigns 
\begin{equation} \label{e:wq}
 \begin{split}
  &\hat{1}= 1, \  \hat{p}_k^i= \sqrt{\hbar} \frac{\p}{\p q^i_k}, 
   \hat{q}^i_k = q^i_k / {\sqrt{\hbar}}, \\
  &(p^i_k p^j_l) \hat{\ } = \hat{p}^i_k \hat{p}^j_l 
	=\hbar \frac{\p}{\p q^i_k} \frac{\p}{\p q^j_l}, \\
   &(p^i_k q^j_l) \hat{\ } = q^j_l \frac{\p}{\p q^i_k},\\ 
  &(q^i_k q^j_l) \hat{\ } = \hat{q}^i_k \hat{q}^j_l /\hbar ,
 \end{split}
\end{equation}

Note that one often has to quantize symplectic instead of infinitesimal
symplectic transformations. Following the common practice in physics, we define
\begin{equation} \label{e:q}
  (e^{A})\hat{}\ := e^{(A) \hat{}\ } ,
\end{equation}
for the symplectic transformation $e^{A(z)}$ in the twisted loop group.

When one restricts the attention to semisimple Frobenius (formal) manifolds,
the situation is even simpler. Let $H_0=\cc^N$ be the Frobenius manifold
such that the orthonormal basis $\{ \phi_{\mu} \}$ form idempotents of the 
product
\[
  \phi_{\mu} * \phi_{\nu} = \delta_{\mu \nu} \phi_{\mu}.
\]
Or in terms of $F_0$
\[
  F_0 (\bt) = \frac{1}{6} \sum_{\mu} (\bt_0^{\mu})^3.
\]
Givental's $\tau$-function \eqref{e:taug} for $H_0$ becomes 
$\tau_G^{H_0} (\bq) = \prod_{\mu} \tau^{pt} (\bq^{\mu})$. It follows from 
Givental's theory that $\tau_G^H$ for any semisimple Frobenius manifold $H$ is
\[
 \tau_G^H = \hat{O}_H \prod \tau_G^{H_0}.
\]
Furthermore, $\hat{O}_H$ is actually an element in the quantized twisted loop 
group. 
\footnote{This statement is not completely correct as the quantized twisted loop
group is not a group. Nevertheless, it is a good heuristic picture.}
Therefore, we conclude that 

\begin{mainlemma}
In order to show that a set of tautological equations holds for $G_g$, 
it suffices to show that this set of tautological equations is invariant 
under quantized loop group action.
\end{mainlemma}

This lemma is our main technical tool to prove Theorem~\ref{t:g1} and 
\ref{t:g2} and therefore Main theorem. In fact, in order to prove the 
invariance of the tautological equations, it is enough to prove the
\emph{infinitesimal} invariance of the tautological equations by the
definition \eqref{e:q}.

\begin{remark}
Morally, one can consider the space of all semisimple Frobenius manifold
as a homogeneous space of quantized twisted loop groups. However, there are
many issues, including the issue of completion alluded before, which make this
assertion invalid.
\end{remark}

\section{Invariance of tautological equations under the action of the twisted
loop groups}

\subsection{Quantization of twisted loop groups}
The twisted loop group is generated by ``lower triangular subgroup'' 
and the ``upper triangular subgroup''.
The lower triangular subgroup consists of $\on{End}(H)$-valued formal formal 
series in $z^{-1}$ $S(z^{-1})= e^{s(z^{-1})}$
satisfying $S^*(-z) S(z) = \mathbf{1}$ or equivalently 
\[
 s^*(-z^{-1}) + s(z^{-1}) =0.
\]
The upper triangular subgroup consists of the regular part of the twisted loop
groups $R(z)= e^{r(z)}$ satisfying $R^*(-z) R(z) = \mathbf{1}$ or equivalently
\begin{equation} \label{e:r}
  r^*(-z) + r(z) =0.
\end{equation}

For illustration, let us work out the quantization of the upper triangular
subgroups. The quantization of $r(z)$ is
\[
 \begin{split}
 \hat{r}(z) = &\sum_{l=1}^{\infty} \sum_{n=0}^{\infty}
    \sum_{i,j} (r_l)_{ij} q^j_n \p_{q^i_{n+l}} \\
  + &\frac{\hbar}{2} \sum_{l=1}^{\infty} \sum_{m=0}^{l-1}
   (-1)^{m+1} \sum_{i j} (r_l)_{ij} \p_{q^i_{l-1-m}} \p_{q^j_m}.
 \end{split}
\]

Let $\frac{d \tau_G}{d \epsilon_r} = \hat{r}(z) \tau_G$. Then

\begin{equation} \label{e:0}
 \begin{split}
  \frac{d}{d \epsilon_r} &\langle \p^{i_1}_{k_1} \p^{i_2}_{k_2} \ldots \rangle
 = \sum_{l=1}^{\infty} \sum_{n=0}^{\infty} \sum_{i,j} (r_l)_{ij} q^j_n
   \langle \p^i_{n+l} \p^{i_1}_{k_1} \ldots \rangle \\
  + &\sum_{l=1}^{\infty} \sum_{i,a} (r_l)_{i i_a} \langle \p^i_{k_a+l}
   \p^{i_1}_{k_1} \ldots \hat{\p^{i_a}_{k_a}} \ldots \rangle \\
  + &\frac{1}{2} \sum_{l=1}^{\infty} \sum_{m=0}^{l-1} (-1)^{m+1}
   \sum_{i j} (r_l)_{ij} \p^{i_1}_{k_1} \p^{i_2}_{k_2} \ldots
   ( \langle \p^i_{l-1-m} \rangle \langle \p^j_m \rangle ) .
 \end{split}
\end{equation}

For $g \ge 1$
\begin{equation} \label{e:g}
 \begin{split}
  &\frac{d \langle \p^{i_1}_{k_1} \p^{i_2}_{k_2} \ldots \rangle_g}
        {d \epsilon_r} \\
  = &\sum_{l=1}^{\infty} \sum_{n=0}^{\infty} \sum_{i,j} (r_l)_{ij} q^j_n
   \langle \p^i_{n+l} \p^{i_1}_{k_1} \ldots \rangle_g \\
  + &\sum_{l=1}^{\infty} \sum_{i,a} (r_l)_{i i_a} \langle \p^i_{k_a+l}
    \p^{i_1}_{k_1} \ldots \hat{\p^{i_a}_{k_a}} \ldots \rangle_g \\
  + &\frac{1}{2} \sum_{l=1}^{\infty} \sum_{m=0}^{l-1} (-1)^{m+1}
   \sum_{i j} (r_l)_{ij} \langle \p^i_{l-1-m} \p^j_m
    \p^{i_1}_{k_1} \p^{i_2}_{k_2} \ldots \rangle_{g-1} \\
  + &\frac{1}{2} \sum_{l=1}^{\infty} \sum_{m=0}^{l-1} (-1)^{m+1}
   \sum_{i j} \sum_{g'=g}^g (r_l)_{ij} \p^{i_1}_{k_1} \p^{i_2}_{k_2} \ldots
   ( \langle \p^i_{l-1-m} \rangle_{g'} \langle \p^j_m \rangle_{g-g'} ).
 \end{split}
\end{equation}

\subsection{Invariance theorems}

\begin{theorem} \emph{($S$-invariance theorem)}
All tautological relations are invariant under action of lower triangular
subgroups of the twisted loop groups.
\end{theorem}

This theorem combined with observations in \cite{GL} and the result from 
\cite{AG2}, which shows that there exists a lower triangular element $S$ 
such that $\hat{S} \tau^X$ is the generating function of ``ancestors'', 
implies that one may assume that $\tau^X$ is the ancestor potential and that 
$q_0=0$ in the proof of $R$-invariance.

Let us use the term ``genus zero relations'' for genus zero dilaton equation,
string equation, and TRR, the term ``genus one relations'' for genus one
Getzler's equation and genus one TRR, the term ``genus two relations''
for genus two equations by Mumford, Getzler, and BP.

\begin{theorem} \emph{($R$-invariance theorem)}
The union of the sets of genus $g'$ relations for $g' \le g$ is invariant 
under the action of upper triangular subgroup, for $g \le 2$.
\end{theorem}

In fact, a stronger ``filtered'' statement holds. We will state the genus two 
part:
\begin{enumerate}
\item The combination of genus zero relations, genus one 
relations and Mumford's equation is $R$-invariant.
\item The combination of genus zero relations, genus one 
relations and genus two Mumford's and Getzler's equations is $R$-invariant.
\end{enumerate}

\begin{remark}
$R$-invariance theorem is expected to hold for all $g$. This will be discussed
in a separate paper.
\end{remark}

\end{document}